\documentclass{amsart}
\begin{document}
\title[Arithmetic analogues]{Arithmetic analogues of some basic concepts from Riemannian geometry}
\author{Alexandru Buium}
\def \Rp{R_p}
\def \Rpi{R_{\pi}}
\def \dpi{\d_{\pi}}
\def \bT{{\bf T}}
\def \cI{{\mathcal I}}
\def \cH{{\mathcal H}}
\def \cJ{{\mathcal J}}
\def \ZN{\bZ[1/N,\zeta_N]}
\def \tA{\tilde{A}}
\def \o{\omega}
\def \tB{\tilde{B}}
\def \tC{\tilde{C}}
\def \alph{A}
\def \bet{B}
\def \bsigma{\bar{\sigma}}
\def \y{^{\infty}}
\def \Ra{\Rightarrow}
\def \uBS{\overline{BS}}
\def \lBS{\underline{BS}}
\def \lB{\underline{B}}
\def \<{\langle}
\def \>{\rangle}
\def \hL{\hat{L}}
\def \cU{\mathcal U}
\def \cF{\mathcal F}
\def \S{\Sigma}
\def \st{\stackrel}
\def \sd{Spec_{\d}\ }
\def \pd{Proj_{\d}\ }
\def \s{\sigma_2}
\def \i{\sigma_1}
\def \bs{\bigskip}
\def \cD{\mathcal D}
\def \cC{\mathcal C}
\def \cT{\mathcal T}
\def \cK{\mathcal K}
\def \cX{\mathcal X}
\def \sX{X_{set}}
\def \cY{\mathcal Y}
\def \cS{X}
\def \cR{\mathcal R}
\def \cE{\mathcal E}
\def \tcE{\tilde{\mathcal E}}
\def \cP{\mathcal P}
\def \cA{\mathcal A}
\def \cV{\mathcal V}
\def \cM{\mathcal M}
\def \cL{\mathcal L}
\def \cN{\mathcal N}
\def \tcM{\tilde{\mathcal M}}
\def \caS{\mathcal S}
\def \cG{\mathcal G}
\def \cB{\mathcal B}
\def \tG{\tilde{G}}
\def \cF{\mathcal F}
\def \h{\hat{\ }}
\def \hp{\hat{\ }}
\def \tS{\tilde{S}}
\def \tP{\tilde{P}}
\def \tA{\tilde{A}}
\def \tX{\tilde{X}}
\def \tcS{\tilde{X}}
\def \tT{\tilde{T}}
\def \tE{\tilde{E}}
\def \tV{\tilde{V}}
\def \tC{\tilde{C}}
\def \tI{\tilde{I}}
\def \tU{\tilde{U}}
\def \tG{\tilde{G}}
\def \tu{\tilde{u}}
\def \chu{\check{u}}
\def \tx{\tilde{x}}
\def \tL{\tilde{L}}
\def \tY{\tilde{Y}}
\def \d{\delta}
\def \e{\chi}
\def \bW{\mathbb W}
\def \bV{{\mathbb V}}
\def \bF{{\bf F}}
\def \bE{{\bf E}}
\def \bC{{\bf C}}
\def \bO{{\bf O}}
\def \bR{{\bf R}}
\def \bA{{\bf A}}
\def \bB{{\bf B}}
\def \cO{\mathcal O}
\def \ra{\rightarrow}
\def \bx{{\bf x}}
\def \f{{\bf f}}
\def \bX{{\bf X}}
\def \bH{{\bf H}}
\def \bS{{\bf S}}
\def \bF{{\bf F}}
\def \bN{{\bf N}}
\def \bK{{\bf K}}
\def \bE{{\bf E}}
\def \bB{{\bf B}}
\def \bQ{{\bf Q}}
\def \bd{{\bf d}}
\def \bY{{\bf Y}}
\def \bU{{\bf U}}
\def \bL{{\bf L}}
\def \bQ{{\bf Q}}
\def \bP{{\bf P}}
\def \bR{{\bf R}}
\def \bC{{\bf C}}
\def \bD{{\bf D}}
\def \bM{{\bf M}}
\def \bZ{{\mathbb Z}}
\def \xtoleqr{x^{(\leq r)}}
\def \hU{\hat{U}}
\def \k{\kappa}
\def \ee{\overline{p^{\k}}}

\newtheorem{THM}{{\!}}[section]
\newtheorem{THMX}{{\!}}
\renewcommand{\theTHMX}{}
\newtheorem{theorem}{Theorem}[section]
\newtheorem{corollary}[theorem]{Corollary}
\newtheorem{lemma}[theorem]{Lemma}
\newtheorem{proposition}[theorem]{Proposition}
\theoremstyle{definition}
\newtheorem{definition}[theorem]{Definition}
\theoremstyle{remark}
\newtheorem{remark}[theorem]{Remark}
\newtheorem{example}[theorem]{\bf Example}
\numberwithin{equation}{section}
\address{University of New Mexico \\ Albuquerque, NM 87131}
\email{buium@math.unm.edu} 
\maketitle


\begin{abstract}
Following recent work of the author, partly in collaboration with T. Dupuy and M. Barrett,
we describe  arithmetic analogues of some key concepts from Riemannian geometry such as: metrics, Chern connections, curvature, etc. Theorems are stated to the effect that the spectrum of the integers has a non-vanishing curvature. 
\end{abstract}

\section{Introduction}

In previous work (initiated in \cite{char} and partly summarized in \cite{book,statupdated})
the author has developed an arithmetic analogue of  differential calculus and, in particular, of differential equations. As explained in \cite{borgerf1, manin}, this theory can be viewed as an alternative approach to ``absolute geometry" (or the ``geometry over the field with one element, ${\mathbb F}_1$") and led to a series of diophantine applications \cite{statupdated}. Once an arithmetic analogue of differential calculus is available one can ask for arithmetic analogues of the basic concepts of differential geometry and, in particular, of Riemannian geometry. Such analogues were recently proposed in \cite{adel1,adel2,adel3,curvature1,curvature2} and led to the somewhat surprising conclusion that the spectrum of the integers, $Spec\ \bZ$,
can be viewed as an (infinite dimensional) manifold which is naturally ``curved" (although, as we shall see, only ``mildly" curved). The aim of this note is to present, in a self contained manner,  some of the ideas and results of this ``arithmetic Riemannian" theory.
For the details of the theory we refer to the papers cited above.

\bigskip

{\it Acknowledgement}.
The author would like to acknowledge partial support from the Simons Foundation
(award 311773)  and from the Romanian National Authority
for Scientific Research (CNCS - UEFISCDI, project number
PN-II-ID-PCE-2012-4-0201).

\section{Main concepts and results}

The best way to present our material is by analogy with  classical differential geometry. In classical differential geometry one starts with an $m$-dimensional smooth manifold $M$ and its ring of smooth functions $C^{\infty}(M)$. 
For our purposes it is enough to think of $M$ as being the Euclidean space $M={\mathbb R}^m$. Also
we would like to think of the  dimension $m$ as going to infinity, $m\ra \infty$. In this paper the arithmetic analogue of ${\mathbb R}^m$, with $m\ra \infty$,  will be the scheme $Spec\ \bZ$; hence the arithmetic analogue of the ring
\begin{equation}
\label{Cinfty}
A:=C^{\infty}({\mathbb R}^m)\end{equation}
 will be the ring of integers $\bZ$ or, more generally,  the ring 
 \begin{equation}
 \label{A}
 A:=\bZ[1/N_0,\zeta_N]\end{equation}
  where $N_0$ is an even integer, $N$ is an integer,  and $\zeta_N$ is a primitive $N$-th root of unity. Let 
  \begin{equation}
  \label{xxx}
  x_1,x_2,...,x_m\end{equation}
   be the coordinates on ${\mathbb R}^m$; then the arithmetic analogues of these coordinates will be a sequence  of primes 
   \begin{equation}
   \label{P}
   {\mathcal P}=\{p_1,p_2,p_3,...\}.\end{equation}
    One can take all primes or, better, all primes not dividing $N_0N$. Next one considers the partial derivative operators 
\begin{equation}
\label{partial}
\frac{\partial}{\partial x_i}:C^{\infty}({\mathbb R}^{\infty})\ra C^{\infty}({\mathbb R}^{\infty}),\ \ i\in \{1,...,m\}.\end{equation}
Following \cite{char} we propose to take, as an analogue of \ref{partial}, the operators 
\begin{equation}
\label{fermat}
\d_{p}:A\ra A,\ \ \d_{p}(a)=\frac{\phi_{p}(a)-a^{p}}{p},\ \ p\in {\mathcal P},
\end{equation}
where $\phi_{p}:A\ra A$ is the unique ring automorphism  of $A$ sending $\zeta_N$ into $\zeta_N^{p}$. More generally the concept of {\it derivation} on a ring $B$ (by which we mean an additive map $B\ra B$ that satisfies the Leibniz rule) has, as an arithmetic analogue, the concept of {\it $p$-derivation} on a ring $B$ (by which we mean a set theoretic map $\d_p:B\ra B$ with the property that the map $\phi_p:B\ra B$ defined by $\phi_p(b)=b^p+p\d_pb$ is a ring homomorphism; we will always denote by $\phi_p$ the ring homomorphism attached to a $p$-derivation $\d_p$ and we shall refer to $\phi_p$ as the {\it  Frobenius lift} attached to $\d_p$). 

The next step in the classical theory is to consider a vector bundle $\pi:E\ra M$ of rank $n$ over the manifold $M$ and to consider the frame bundle $F(E)$ of  $E$; a point of $F(E)$ is a point $P$ of $M$ together with a basis of $\pi^{-1}(P)$.  The frame bundle $F(E)$ is a principal homogeneous space for the group $GL_n$; and if $E$ is a trivial vector bundle (which we shall assume from now on) then $F(E)$ is identified with $M\times GL_n$. (Note that the rank $n$ of $E$ and the dimension $m$ of $M$ in this picture are unrelated.) We want to review the classical concept of connection in $F(E)$; we shall do it in a somewhat non-standard way so that the transition to arithmetic becomes more transparent. Indeed consider 
an $n \times n$ matrix $x=(x_{ij})$ of indeterminates and consider
the ring of polynomials over $A$, with the determinant inverted,
\begin{equation}
\label{B}
B=A[x,\det(x)^{-1}].
\end{equation}
Note that $B$ is naturally a subring of the ring $C^{\infty}(M\times GL_n)$. Then by a  {\it connection} on $F(E)=M\times GL_n$ we will understand a tuple $\d=(\d_i)$ of derivations 
\begin{equation}
\label{deltai}
\d_i:B\ra B,\ \ i\in \{1,...,m\},
\end{equation}
lifting the derivations \ref{partial}. We say the  connection is {\it linear} if $\d_ix=A_ix$ for some $n\times n$ matrix $A_i$ with coefficients in $A$; linearity corresponds to the concept of invariance of classical connections under the right action of $GL_n$ on the frame bundle $F(E)$. For a linear connection $\d$ as above one can define the {\it curvature} as the matrix $(\varphi_{ij})$ of commutators 
\begin{equation}
\label{varphi}
\varphi_{ij}:=[\d_i,\d_j]:B\ra B,\ \ \ i,j\in \{1,...,m\}.\end{equation}
One has then $\varphi_{ij}(x)=F_{ij}x$ where $F_{ij}$ is the matrix given by the classical formula
\begin{equation}
\label{formula}
F_{ij}:=\d_iA_j-\d_jA_i-[A_i,A_j].
\end{equation}
We would like to introduce now an arithmetic analogue of connection and curvature. The first step is clear: we consider a ring $B$ defined as in \ref{B} but where $A$ is given now by \ref{A} rather than by \ref{Cinfty}. A first attempt to define arithmetic connections would be to consider families of $p$-derivations $\d_p:B\ra B$, $p\in {\mathcal P}$, lifting the $p$-derivations \ref{fermat}; one would then proceed by considering their commutators on $B$ (or, if necessary, expressions derived from these commutators). But the point is that the examples of ``arithmetic connections" we will encounter in practice (when we develop arithmetic analogues of the Chern connections of classical differential geometry) will never lead to $p$-derivations $B\ra B$! What we shall be led to is, rather, a concept we next introduce. 
For each $p\in {\mathcal P}$ we consider the $p$-adic completion of $B$:
\begin{equation}
\label{completion}
B^{\widehat{p}}:=\lim_{\leftarrow} B/p^nB.
\end{equation}
Then we define an {\it arithmetic connection} on $GL_n$ to be a family $(\d_p)$ of $p$-derivations
\begin{equation}
\label{arithmetic}
\d_p:B^{\widehat{p}}\ra B^{\widehat{p}},\ \ p\in {\mathcal P},
\end{equation}
lifting the $p$-derivations in \ref{fermat}. We do not impose any condition analogous to linearity; instead, what happens is that our arithmetic connections of interest turn out to enjoy a certain  invariance property with respect to right translations by the elements of the normalizer of the maximal (diagonal) torus of $GL_n$. This invariance can be viewed as a substitute for linearity and will not be discussed here further.
Leaving the linearity issue aside we are facing, at this point,  a more severe dilemma: our $p$-derivations $\d_p$ in \ref{arithmetic} do not act on the same ring, so there is no a priori way of considering their commutators and, hence, it does not seem possible to define, in this way, the notion of curvature.
It will turn out, however, that our arithmetic connections of interest will satisfy an interesting property which we call ``being global along the identity," and which will allow us to define curvature via commutators. Here is the definition of this property. Consider the matrix $T=x-1$, where $1$ is the identity matrix.  We say that an arithmetic connection $(\d_p)$ on $GL_n$, with attached family of Frobenius lifts $(\phi_p)$, is {\it global along $1$} if, for all $p$, $\phi_p$ sends the ideal of $1$ into itself and, moreover,  the induced homomorphism $\phi_p:A^{\widehat{p}}[[T]]\ra A^{\widehat{p}}[[T]]$
sends the ring $A[[T]]$ into itself. If the above holds then one can consider the commutator $[\phi_p,\phi_{p'}]:A[[T]]\ra A[[T]]$ for all $p,p'\in {\mathcal P}$; this commutator is divisible by $pp'$ and one can define the {\it curvature} of $(\d_p)$ as the matrix $(\varphi_{pp'})$ with entries
\begin{equation}
\label{definecurv}
\varphi_{pp'}:=\frac{1}{pp'}[\phi_p,\phi_{p'}]:A[[T]]\ra A[[T]],\ \ \ p,p'\in {\mathcal P}.
\end{equation}
The idea of comparing $p$-adic phenomena for different $p$'s by ``moving along the identity section" is borrowed from \cite{laplace} where it was referred to as ``analytic continuation along primes."
Of course, in order for the above definitions to be interesting, we will need to: 

1) find natural ``metric" arithmetic connections on $GL_n$, 

2) show that these connections are global along $1$, and

3) compute and interpret the curvatures of these connections.

\noindent We embark now on explaining how this program can be achieved. First we go back to classical differential geometry and we ``recall" the definition of the Chern connection \cite{kobayashi}. We shall present this definition in the ``real setting" only, where the Chern connections should be more appropriately referred to as Duistermaat connections \cite{dui}; for the ``complex setting" we refer to \cite{adel2,curvature1}. So let us consider the ring $A=C^{\infty}({\mathbb R}^m)$  and let $q$ be an $n\times n$ invertible matrix with coefficients in $A$ which is either symmetric ($q^t=q$) or antisymmetric ($q^t=-q$), where the $t$ superscript means {\it transposition}. Of course, a symmetric $q$ as above is viewed as a ``metric" while an antisymmetric $q$ is viewed as a ``$2$-form." Set $G=GL_n$ and consider the  maps of schemes over $A$,
$\cH_q:G\ra G$, $\cB_q:G\times G\ra G$ defined by
$\cH_q(x)=x^tqx$ and $\cB_q(x,y)=x^tqy$. We continue to denote by the same letters the corresponding maps of rings $B\ra B$ and $B\ra B\otimes_A B$.
Consider the {\it trivial} (linear) connection $\d_0=(\d_{0i})$ on $G$ defined by $\d_{0i}x=0$. Then one can easily check (see below)  that there is a unique 
 linear connection $(\d_i)$ on $G$ such that the following diagrams are commutative:
 \begin{equation}
 \label{got}
 \begin{array}{rcl}
 B & \stackrel{\d_i}{\longleftarrow} & B\\
 \cH_q \uparrow &\ &\uparrow \cH_q\\
 B & \stackrel{\d_{0i}}{\longleftarrow} & B\end{array}
 \ \ \ \ \ \ \ 
  \begin{array}{rcl}
 B & \stackrel{\d_i\otimes 1+1\otimes \d_{0i}}{\longleftarrow} & B\otimes_A B\\
 \d_{0i}\otimes 1+1\otimes \d_i \uparrow & \  & \uparrow \cB_q\\
 B\otimes_A B & \stackrel{\cB_q}{\longleftarrow} & B\end{array}
 \end{equation}
 This $\d$ can be referred to as the {\it Chern connection} attached to $q$.
 The definition just given may look  non-standard. It turns out that the Chern connection we just defined is a real analogue \cite{dui} of the usual Chern connection in differential geometry \cite{kobayashi} (in which $\d_0$ is an analogue of a complex structure). To see this set
 $\Gamma_i = -A_i^t$, let $\Gamma_{ij}^k$ be the $(j,k)$-entry of $\Gamma_i$
(the Cristoffel symbols), and set $\Gamma_{ijk} :=\Gamma_{ij}^lq_{lk}$ (Einstein notation). Assume we are in the symmetric case, $q^t=q$. Then 
the commutativity of the left diagram in \ref{got} is equivalent to the condition
\begin{equation}
\label{horizontal}
\d_i q_{jk}=\Gamma_{ijk}+ \Gamma_{ikj},
\end{equation}
and the commutativity of the right diagram in \ref{got} is equivalent to the condition
\begin{equation}
\label{symmetric}
\Gamma_{ijk}=\Gamma_{ikj};
\end{equation}
so the Chern connection attached to $q$ is given by
\begin{equation}
\label{solution}
\Gamma_{ijk}=\frac{1}{2}\d_i q_{jk}.
\end{equation}
The Chern connection will have an arithmetic analogue to be explained presently.
The condition \ref{horizontal} expresses the fact that $q$ is {\it parallel} with respect to the connection $\d$.
It is important to note, however,  that, in our setting, the {\it torsion} is not defined and, in particular, the symmetry in \ref{symmetric} has nothing to do with the vanishing of the torsion.  On the other hand, if one takes $E$ to be the tangent bundle of $M$ (so in particular $n=m$),  then the condition that the {\it torsion of $\d$ vanishes}  is given by:
\begin{equation}
\label{torsion}
\Gamma_{ijk}=\Gamma_{jik}
\end{equation}
which is a symmetry condition rather different from \ref{symmetric}. By the way there is a unique connection $\d$ such that conditions \ref{horizontal} and \ref{torsion} are satisfied; this connection is referred to as the {\it Levi-Civita connection} and is given by the formula
\begin{equation}\label{windy}
\Gamma_{kij}=\frac{1}{2}\left(\d_k q_{ij}+\d_iq_{jk}-\d_jq_{ki}\right).
\end{equation}
The Levi-Civita connection does not seem to have an arithmetic analogue in our theory. 

Now we move to the arithmetic situation. So let $A=\bZ[1/N_0,\zeta_N]$.
Let $q\in GL_n(A)$ with $q^t=\pm q$. Set $G=GL_n=Spec\ B$, viewed as a group scheme over $A$. Attached to $q$ we have, again,   maps $\cH_q:G\ra G$ and  $\cB_q:G\times G\ra G$ defined by
$\cH_q(x)=x^tqx$ and $\cB_q(x,y)=x^tqy$. 
We continue to denote by $\cH_q,\cB_q$ the maps induced on the $p$-adic completions $G^{\widehat{p}}$ and $G^{\widehat{p}}\times G^{\widehat{p}}$. Consider  the unique arithmetic connection $\d_0=(\d_{0,p})$ on $G$ with $\d_{0,p}x=0$ and denote by $(\phi_p)$ and $(\phi_{0,p})$ the families of lifts of Frobenius attached to $\d$ and $\d_0$ respectively. Then one has the following:

\begin{theorem}
\label{existence}
\cite{adel2}.
For any $q\in GL_n(A)$ with $q^t=\pm q$ there exists a unique arithmetic connection $\d$ such that the following diagrams are commutative:
$$
  \begin{array}{rcl}
G^{\widehat{p}} & \stackrel{\phi_p}{\longrightarrow} & G^{\widehat{p}}\\
\cH_q  \downarrow & \  & \downarrow \cH_q \\
G^{\widehat{p}} & \stackrel{\phi_{0,p}}{\longrightarrow} & G^{\widehat{p}}\\
\end{array}\ \ \ \ \ \ 
\begin{array}{rcl}
G^{\widehat{p}} & \stackrel{\phi_{0,p} \times \phi_p}{\longrightarrow} & G^{\widehat{p}}\times G^{\widehat{p}}\\
\phi_p \times \phi_{0,p} \downarrow & \  & \downarrow \cB_q\\
 G^{\widehat{p}}\times G^{\widehat{p}} & \stackrel{\cB_q}{\longrightarrow} & G^{\widehat{p}}\end{array}$$
 \end{theorem}
 
The arithmetic connection $\d$ is referred to as the  {\it Chern connection} (on $G=GL_n$) attached to $q$. 
The various $q$'s with $q^t=\pm q$ lead to the various forms of the classical groups $Sp_n$ and $SO_n$. A similar theorem is proved in \cite{adel2} for the classical groups $SL_n$. 
Note the following relation between the ``Christoffel symbols" defining our Chern connection and the Legendre symbol. We explain this in a special case.
 Let $q\in GL_1(A)=A^{\times}$, $A=\bZ[1/N_0]$, and let $\d=(\d_p)$ be the Chern connection  associated to $q$.
   Then it turns out that $\phi_p:G^{\widehat{p}}\ra G^{\widehat{p}}$ is defined by $\phi_p:\bZ_p[x,x^{-1}]^{\widehat{p}}\ra \bZ_p[x,x^{-1}]^{\widehat{p}}$,
   \begin{equation}
   \phi_p(x)= 
   q^{(p-1)/2}\left(\frac{q}{p}\right)x^p,\end{equation}
where $\left(\frac{q}{p}\right)$ is the Legendre symbol of $q\in A^{\times}\subset \bZ_{(p)}$.

Next one can ask which of these arithmetic connections admit curvatures.
One  has:

\begin{theorem}\cite{curvature1}.
 If all the entries of $q$ are roots of unity or $0$ then
 the Chern connection $\d$ attached to $q$ is global along $1$. In particular $\d$ has a well defined curvature.\end{theorem}
 
 Next we address the question of computing the curvature of Chern connections. Let us say that a matrix $q\in GL_n(A)$ is  split if it is one of the following:
 \begin{equation}
 \label{scorpion3}
  \left(\begin{array}{cl} 0 & 1_r\\-1_r & 0\end{array}\right),\ \ 
\left( 
\begin{array}{ll} 0 & 1_r\\1_r & 0\end{array}\right),\ \ 
\left( \begin{array}{lll} 1 & 0 & 0\\
0 & 0 & 1_r\\
0 & 1_r & 0\end{array}\right),
\end{equation}
where $1_r$ is the $r\times r$ identity matrix and $n=2r,2r,2r+1$ respectively.
These are matrices that define the classical split groups $Sp_{2r},SO_{2r},SO_{2r+1}$, respectively. 
One has the following:

 \begin{theorem}
 \label{coconut}
 \cite{curvature1}.
 Let $q$ be  split  and let $(\varphi_{pp^{\prime}})$ be the curvature of the Chern connection on $G$ attached to $q$.
 
 1) Assume  $n\geq 4$. Then for all $p\neq p^{\prime}$ we have $\varphi_{pp^{\prime}}\neq 0$.
 
 2) Assume $n$ even. Then for all $p,p^{\prime}$ we have $\varphi_{pp^{\prime}}(T)\equiv 0$ mod $(T)^3$.
 
 3) Assume $n=2$ and $q^t=-q$. Then for all $p,p^{\prime}$ we have $\varphi_{pp^{\prime}}= 0$.
 
 4) Assume $n=1$. Then for all $p,p'$ we have $\varphi_{pp'}=0$.
  \end{theorem}
  
In assertion 2) we denoted by $(T)^3$  the cube of the ideal in $A[[T]]$ generated by the entries of the matrix $T$. Assertion 1) morally says that $Spec\ \bZ$ is ``curved," while assertion
2) morally says that $Spec\ \bZ$ is only ``mildly curved." Note that the theorem says nothing about  the vanishing of the  curvature $\varphi_{pp'}$ in case $n=2,3$ and  $q^t=q$; our method of proof does not seem to apply to these cases.

The theory explained above has a ``complex analogue" (or, rather, a $(1,1)$-analogue) for which we refer to \cite{adel2,curvature1}. 
This theory (that largely follows \cite{curvature1})  was based on what we called ``analytic continuation between primes"; this was the key to making Frobenius lifts corresponding to different primes act on a same ring.  There is a different approach towards making Frobenius lifts comparable; this approach  was developed in \cite{curvature2}.
The idea in \cite{curvature2} was to show that if $\d=(\d_p)$ is the Chern connection attached to a matrix $q\in GL_n(A)$ with $q^t=\pm q$, then one can find 
 correspondences $\Gamma_p=(Y_p,\pi_p,\varphi_p)$ on $G=GL_n$, i.e.  maps of $A$-schemes
\begin{equation}
\label{corr}
\pi_p:Y_p\ra G,\ \ \ \varphi_p:Y_p\ra G,
\end{equation}
 such that:
 
 i) $\pi_p$ are affine and \'{e}tale, 
 
 ii) $\pi^{\widehat{p}}_p:Y_p^{\widehat{p}}\ra G^{\widehat{p}}$ are isomorphisms, and 
 
 iii) $\varphi_p^{\widehat{p}}=\phi_p\circ \pi_p^{\widehat{p}}:Y_p^{\widehat{p}}\ra G^{\widehat{p}}$.
 
\noindent In other words the correspondences \ref{corr} are ``algebraizations" 
of our Frobenius lifts $\phi_p$; the system $(\Gamma_p)$ is referred to as a {\it correspondence structure} for $(\d_p)$; it is not unique but does have some ``uniqueness features" (cf. \cite{curvature2}).  On the other hand any correspondence $\Gamma_p$ as in \ref{corr}  acts on the field $E$ of rational functions of $G$ by the formula $\Gamma_p^*:E\ra E$, 
\begin{equation}
\Gamma_p^*(z)=\text{Tr}_{\pi_p}(\varphi_p^*(z)),\ \ z\in E,
\end{equation}
where $\text{Tr}_{\pi_p}:F_p\ra E$ is the trace of the extension $\pi_p^*:E\ra F_p:=Y_p\otimes_G E$ and $\varphi_p^*:E\ra F_p$ is induced by $\varphi_p$.
By the way the degrees of the extensions $\pi_p^*:E\ra F_p$ and $\varphi_p^*:E\ra F_p$ will be referred to as the {\it left degree} and the {\it right degree} of $\Gamma_p$ respectively. Also we say $\Gamma_p$ is {\it irreducible} if $F_p$ is a field.
 So one can define the {\it $*$-curvature} of the arithmetic connection $(\d_p)$ as the matrix $(\varphi^*_{pp'})$ where 
\begin{equation}
\varphi^*_{pp'}:=\frac{1}{pp'}[\Gamma^*_p,\Gamma^*_{p'}]:E\ra E,\ \ \ p,p'\in {\mathcal P}.
\end{equation}
Note that, in this way, we have defined a concept of ``curvature" for Chern connections attached to arbitrary $q$'s (that  do not necessarily have entries zeroes or roots of unity). 
There is a $(1,1)$-version of the above as follows.
 Given one more arithmetic connection $\overline{\d}=(\overline{\d}_p)=:(\d_{\overline{p}})$ with correspondence structure $(\overline{\Gamma}_p)=:(\Gamma_{\overline{p}})$ one can define the {\it $(1,1)$-$*$-curvature} of $(\Gamma_p)$ with respect to $(\Gamma_{\overline{p}})$ as the family $(\varphi_{p\overline{p}'}^*)$ where $\varphi^*_{p\overline{p}'}$ is the additive endomorphism
  \begin{equation}
  \label{varfistar2}
  \varphi^*_{p\overline{p}'}:=\frac{1}{pp'}[\Gamma^*_{\overline{p}'},\Gamma^*_p]:E\ra E\ \ \text{for $p\neq p'$, and}
 \ \ \ 
  \varphi^*_{p\overline{p}}:=\frac{1}{p}[\Gamma^*_{\overline{p}},\Gamma^*_p]:E\ra E.
  \end{equation}
 In what follows we  let $\overline{\d}$ be equal to $\d_0=(\d_{0,p})$, where $\d_{0,p}x=0$; we give $\overline{\d}$ the correspondence structure 
$(\Gamma_{\overline{p}})=(G,\pi_{\overline{p}},\varphi_{\overline{p}})$, $\pi_{\overline{p}}$ the identity, and $\varphi_{\overline{p}}(x)=x^{(p)}$.

 \begin{theorem}
 \label{alta}
 \cite{curvature2}.
 
    1) Assume $n=2$ and 
   $q$ is split with $q^t=-q$.
   Then $\Gamma_p$ is irreducible and has  left degree $2$ and right degree $2p^4$. Moreover the $*$-curvature satisfies $\varphi^*_{pp'}=0$ for all $p,p'$ while the $(1,1)$-$*$-curvature satisfies $\varphi^*_{p\overline{p}'}\neq 0$ for all $p,p'$.
   
   2)  Assume $n=2$ and 
   $q$ is split with $q^t=q$.
  Then $\Gamma_p$ is irreducible and  has left degree $4$. Moreover the $(1,1)$-$*$-curvature satisfies $\varphi_{p\overline{p}'}^*\neq 0$ for all $p,p'$. 
 \end{theorem}
 
Once again, the theorem says nothing about the $*$-curvature in case $n=2$ and $q^t=q$; our method of proof  does not seem to apply to this case.

\section{Final remarks}

The theory outlined above  should be viewed as a first step in a  program of developing a {\it differential geometry of $Spec\ \bZ$}. Other types of curvature ({\it Ricci, mean, scalar}) are developed in \cite{curvature1} and lead to some interesting Dirichlet series. 
An {\it arithmetic Maurer-Cartan connection} and a Galois theory attached to it is given in \cite{adel2,adel3}; this Galois theory should be viewed as an {\it arithmetic gauge theory} and should be further developed. It might be possible to attach 
{\it deRham cohomology classes}  to our curvatures and to link them  to the \'{e}tale cohomology of $Spec\ \bZ$. Links between arithmetic connections and  Galois representations might exist that mimic the link between flat connections on vector bundles over manifolds and representations of the fundamental group of those manifolds. We hope to come back to these issues in future work.

\end{document}